\documentclass[11pt]{amsart}
\usepackage{graphicx}
\newtheorem{thm}{Theorem}[section]

\newtheorem{lem}[thm]{Lemma}

\theoremstyle{definition}

\theoremstyle{remark}

\begin{document}

\title[Time-inhomogeneous birth-and-death processes]{Conditions for recurrence and transience for time-inhomogeneous birth-and-death processes} %
\author[Abramov]{Vyacheslav M. Abramov$^*$}
\address{24 Sagan Drive, Cranbourne North, Victoria-3977, Australia}%

\email{vabramov126@gmail.com}%

\thanks{$^*$ https://orcid.org/0000-0002-9859-100X}%
%\date{}%
%\dedicatory{}%
%\commby{}%

\begin{abstract}
We derive the conditions for recurrence and transience for time-inhomogeneous birth-and-death processes considered as random walks with positively biased drifts. We establish a general result, from which the earlier known particular results by Menshikov and Volkov [\emph{Electron. J. Probab.}, \textbf{13} (2008), paper No. 31, 944--960] follow.
\end{abstract}

\maketitle

\section{Introduction}\label{S0}

Menshikov and Volkov \cite{MV} studied the time-inhomogeneous random walk $X_t\in\mathbb{R}$, $t=1,2,\ldots$, with the drift satisfying the condition
\begin{equation}\label{1}
\mathsf{E}\{X_{t+1}-X_t~|~X_t=x\}\sim \rho \frac{|x|^\alpha}{t^\beta},
\end{equation}
where $\rho>0$, $\alpha$ and $\beta$ were some constants, and the meaning of the symbol `$\sim$' was further clarified for the case studies of the paper. There were some assumptions (indicated below) and the conditions on the parameters $\alpha$, $\beta$ and $\rho$ in \cite{MV}, under which the random walk was recurrent or transient. The study of this random walk was associated with urn models and random walks with vanishing drifts that have been intensively studied in the literature (see the references in \cite{MV}), and the obtained results covered many known models from these areas.

The basic assumptions, under which the study in \cite{MV} has been provided for all considered cases, are quite natural: (H1) uniform boundedness of jumps, (H2) uniform non-degeneracy on $[a,\infty)$, and (H3) uniform boundedness of time to leave $[0,a]$ (for further details see \cite{MV}).

In the case when $\alpha<\beta$, it was shown that the random walk is transient, if $\{(\alpha,\beta): 0\leq\beta < 1, 2\beta-1 < \alpha < \beta\}$, and it is recurrent otherwise. In the case when $\alpha=\beta=1$ the random walk started at zero was described by the equalities
\begin{eqnarray}\label{2}
\mathsf{P}\{X_{t+1}=n\pm1~|~X_t=n\}&=&\frac{1}{2}\pm\rho\frac{n}{2t}, \quad  1\leq n\leq t,\\
\mathsf{P}\{X_{t+1}=1~|~X_t=0\}&=&1.  \nonumber
\end{eqnarray}
It was shown that it is transient for $\rho>1/2$ and recurrent for $\rho<1/2$.

Along with the indicated cases, \cite{MV} also discussed a number of different cases, one of which was left unstudied, since the martingale methodology used there seems to be insufficient for the resolving that case.  The boundary configurations of the mentioned case are $\alpha=2\beta-1$ and $\beta\in\left(0,\frac{1}{2}\right)\cup\left(\frac{1}{2},1\right)$. This motivates to find a new approach to this problem to cover this case, in particular. With that new approach, we shall study a generalized version of \eqref{2}:
\begin{eqnarray}\label{4}
\mathsf{P}\{X_{t+1}=n\pm1~|~X_t=n\}&=&\frac{1}{2}\pm\varphi(n,t), \quad  1\leq n\leq t,\\
\mathsf{P}\{X_{t+1}=1~|~X_t=0\}&=&1,\nonumber
\end{eqnarray}
for which the required assumptions will be provided later in the paper.

The random walks defined by \eqref{1} and \eqref{2} have positively biased drift. Considerations provided in \cite{MV} implicitly assumed that the drifts considered in the random walks there all are vanishing as $t\to\infty$,
and the conditions for the parameters under which it is true have not been discussed. In the analysis provided in the present paper for some particular examples, this issue is considered as well.

The other possible model of random walks in $\mathbb{R}$ with positively biased drift when a random walk takes positive values and negatively biased drifts when a random walk takes negative values, that is closely related to that considered in \cite{MV} is
\[
\mathsf{E}\{X_{t+1}-X_t~|~X_t=x\}\sim \rho\frac{\mathrm{sign}(x) |x|^\alpha}{t^\beta}, \quad t=1, 2,\ldots,
\]
where the meaning of `$\sim$' is the same as in \cite{MV}.
However the study of this model is similar to that of \eqref{1} because of the symmetry.

Since the processes described by \eqref{1} and \eqref{2} are with positively biased drifts, they can be assumed to be given in $\mathbb{R}_+$. Furthermore, the assumptions (H1), (H2) and (H3) enable us to further specify the random walks assuming that they are time-inhomogeneous processes of the birth-and-death type (or simply time-inhomogeneous birth-and-death processes). That is, we will further assume that the jumps of the random walks take the values $\pm1$ with probabilities depending on state and time, and the continuous time processes take the values in $\mathbb{Z}_+$. This will enable us to approach the model described by \eqref{4}.

For the study of recurrence and transience this assumption is not restrictive. If a process of the birth-and-death type is recurrent (resp. transient), then there is a wide class of closely related Markov chain models satisfying the same property of recurrence (resp. transience), the jumps of which take the values in $\mathbb{Z}$ rather than in $\{-1, 1\}$.

For simplicity of explanations, consider a time-homogeneous (that is ordinary) birth-and-death process with positive birth rates $\lambda_n$ and positive death rates $\mu_n$ ($\lambda_n>\mu_n$, $\lambda_n+\mu_n=1$). Observe this continuous process at the discrete times $t=0, 1, 2,\ldots$. At the initial time $t=0$ the birth-and-death process is in state $0$. Up to time $t=1$, a random number of births and deaths can occur, and the state of the system at $t=1$ is defined as the difference between the numbers of births and deaths before time $t=1$ to be denoted by $n_1$. Similarly, at time $t=2$ the state of the system is denoted $n_2$, etc.

Apparently, for a time-inhomogeneous birth-and-death process (that will be defined later in the paper) the construction of a Markov process with the jumps belonging to $\mathbb{Z}$ is the same. Specifically, we obtain a new Markov chain $Y_t$, the jumps of which belong to $\mathbb{Z}$. For this newly defined Markov chain, \eqref{4} is rewritten as
\[
\mathsf{E}\{Y_{t+1}-Y_t~|~Y_t=n\}\asymp2\varphi(n,t), \quad n\to\infty.
\]

Condition (H1) in \cite{MV} is restrictive. According to (H1) the jumps of a Markov chain are uniformly bounded for all $t$. For the model that is built above, the assumption is weaker. The jumps are not uniformly bounded, but have all moments. The only difference is that the jumps of Markov chains studied in \cite{MV} take the values in $\mathbb{R}$, while the example below suggests the jumps belonging to $\mathbb{Z}$. We think, however, that it is not crucial.

Thus, the aforementioned arguments enable us to conclude that the model described by \eqref{4} is not less general than the model described by
\[
\mathsf{E}\{X_{t+1}-X_t~|~X_t=n\}=\varphi(n,t), \quad t=0, 1, 2,\ldots,
\]
under Condition (H1).

In the present paper, we suggest a new approach to the problem on recurrence and transience for time-inhomogeneous Markov chains. By using the techniques of stochastic calculus, we provide a simpler and transparent proofs of the results under general assumptions on the increments. We first derive a stochastic variant of Chapman--Kolmogorov equations. Then we apply martingale techniques for continuous time processes and reduce the problem to the known criteria of recurrence or transience for birth-and-death processes.
The stochastic variant of the Chapman--Kolmogorov equations appeared for the first time in the paper of Kogan and Liptser \cite{KL}, and then it has been used in many papers (e.g. \cite{A2, A6, A3, KLS}), the majority of which are from the area of queueing networks.

The rest of the paper is structured into four sections. In Section \ref{S1} we formulate the theorem, and in Section \ref{S2} we prove it. The proof of the theorem in Section \ref{S2} consists of three parts included into Sections \ref{S2.1}--\ref{S2.3}. In Section \ref{S2.1}, for a time-inhomogeneous birth-and-death process we derive a stochastic variant of the Chapman--Kolmogorov equations and prove the convergence of the characteristics of the time-inhomogeneous birth-and-death process to the ordinary one. In Section \ref{S2.2} we recall some facts on recurrence and transience of birth-and-death processes known from the literature. In Section \ref{S2.3}, we finalize the proof of the theorem. In Section \ref{S3}, we provide a few examples for the random walks studied in the paper. Some of the examples support the results obtained in \cite{MV}. In Section \ref{S4}, we conclude the paper.

\section{Formulation of the main result}\label{S1}
Let $X_t$ be the random walk defined as
\begin{eqnarray}\label{3}
\mathsf{P}\{X_{t+1}=n\pm1~|~X_t=n\}&=&\frac{1}{2}\pm\varphi(n,t), \quad  1\leq n\leq t.\\
\mathsf{P}\{X_{t+1}=1~|~X_t=0\}&=&1,  \nonumber
\end{eqnarray}
for $t=0, 1,\ldots$. Assume that $\varphi(n,t)$ well-defines the right-hand side of \eqref{3} (i.e. keeps the probability distributions correctly defined), and it is a decreasing in $t$ Borel function.

\begin{thm}\label{T1}
The random walk $X_t$ is recurrent if there exist $c<1$ and $n_0$ such that $\varphi(n,n^2)\leq c/(4n)$ for all $n\geq n_0$, and $X_t$ is transient if there exist $c>1$ and $n_0$ such that $\varphi(n,n^2)\geq c/(4n)$ for all $n\geq n_0$.
\end{thm}

The proof of Theorem \ref{T1} is given in the next section.

\section{Proof of the theorem}\label{S2}
\subsection{Time-inhomogeneous birth-and-death processes}\label{S2.1}
We consider a time-inhomogeneous version of the birth-and-death process $Z(\tau)$ with the instantaneous positive rates $\lambda_{n,\tau}$ and $\mu_{n,\tau}$ in state $n$ at time $\tau$. To emphasize the continuous time process, we use $\tau$ as a continuous time parameter, rather then $t$ that was a discrete time. The process $Z(\tau)$ is assumed to be right-continuous with left limits. The functions $\lambda_{n,\tau}$ and $\mu_{n,\tau}$ are assumed to be continuous in variable $\tau$ with the finite positive limits
\begin{eqnarray}
  \lambda_n &=& \lim_{\tau\to\infty}\lambda_{n,\tau},\label{5} \\
  \mu_n &=& \lim_{\tau\to\infty}\mu_{n,\tau}.\label{7}
\end{eqnarray}
As well, we assume that the existence of the limits below satisfying
\begin{equation}\label{14}
\lim_{n\to\infty}\lambda_n>0, \quad \lim_{n\to\infty}\mu_n>0.
\end{equation}

The system of equations for the process $Z(\tau)$ in terms of the indicators $\mathsf{I}\{Z(\tau)=n\}$  is a stochastic analogue of the Chapman--Kolmogorov equations:
\begin{equation}\label{8}
\begin{aligned}
\mathsf{I}\{Z(\tau)=n\}=&\mathsf{I}\{Z(\tau-)=n+1\}\mathsf{I}\{\mathrm{M}_{n+1}(\tau)-\mathrm{M}_{n+1}(\tau-)=1\}\\
&+\mathsf{I}\{Z(\tau-)=n-1\}\mathsf{I}\{\Lambda_{n-1}(\tau)-\Lambda_{n-1}(\tau-)=1\}\\
&+\mathsf{I}\{Z(\tau-)=n\}\mathsf{I}\{\mathrm{M}_{n,}(\tau)-\mathrm{M}_{n}(\tau-)=0\}\\
& \ \ \ \ \times\mathsf{I}\{\Lambda_{n}(\tau)-\Lambda_{n}(\tau-)=0\},\\
&n=1,2,\ldots
\end{aligned}
\end{equation}
and
\begin{equation}\label{9}
\begin{aligned}
\mathsf{I}\{Z(\tau)=0\}=&\mathsf{I}\{Z(\tau-)=1\}\mathsf{I}\{\mathrm{M}_{1}(\tau)-\mathrm{M}_{1}(\tau-)=1\}\\
&+\mathsf{I}\{Z(\tau-)=0\}\mathsf{I}\{\Lambda_{0}(\tau)-\Lambda_{0}(\tau-)=0\}.
\end{aligned}
\end{equation}
Here $Z(\tau-)$ denotes the value of the process $Z(\tau)$ immediately before the point $\tau$ (if $\tau$ is a point of continuity of $Z(\tau)$, then $Z(\tau-)=Z(\tau)$, otherwise $|Z(\tau)-Z(\tau-)|=1$),
$\Lambda_{n}(\tau)$ and $\mathrm{M}_{n}(\tau)$ for each fixed $n$ denote the time-inhomogeneous Poisson processes with the instantaneous rates $\lambda_{n,\tau}$ and $\mu_{n,\tau}$, respectively.

According to the Doob-Meyer semimartingale decomposition (e.g. \cite{J-Sh, L-Sh, Pr}), written here in the form of stochastic differentials, we have $\mathrm{d}\Lambda_{n}(\tau)=\lambda_{n,\tau}\mathrm{d}\tau+\mathrm{d}\boldsymbol{M}_{\Lambda_{n}(\tau)}$ and, respectively, $\mathrm{d}\mathrm{M}_{n}(\tau)=\mu_{n,\tau}\mathrm{d}\tau+\mathrm{d}\boldsymbol{M}_{\mathrm{M}_{n}(\tau)}$. Here $\lambda_{n,\tau}\mathrm{d}\tau$ and $\mu_{n,\tau}\mathrm{d}\tau$ are the compensators of $\mathrm{d}\Lambda_{n}(\tau)$ and $\mathrm{d}\mathrm{M}_{n}(\tau)$, and $\mathrm{d}\boldsymbol{M}_{\Lambda_{n}(\tau)}$ and $\mathrm{d}\boldsymbol{M}_{\mathrm{M}_{n}(\tau)}$ are the square integrable martingales, respectively, all written in the form of stochastic differentials.

It is not difficult to explain that the compensators for the stochastic differentials $\mathrm{d}\Lambda_{n}(\tau)$ and $\mathrm{d}\mathrm{M}_{n}(\tau)$ have the forms $\lambda_{n,\tau}\mathrm{d}\tau$ and $\mu_{n,\tau}\mathrm{d}\tau$. The compensator of an ordinary Poisson process with rate $\lambda$ is known to be equal to $\lambda\tau$ (e.g. \cite{J-Sh, L-Sh}). That is, for its stochastic differential we have $\lambda\mathrm{d}\tau$. In the case of a time-inhomogeneous Poisson process, the instantaneous rate coincides with the rate of an ordinary Poisson process given in point $\tau$.

Using the aforementioned facts, our system of equations is rewritten as follows:
\[
\begin{aligned}
\mathrm{d}\mathsf{I}\{Z(\tau)=n\}=&\mathsf{I}\{Z(\tau-)=n+1\}\mu_{n+1,\tau}\mathrm{d}\tau+\mathsf{I}\{Z(\tau-)=n-1\}\lambda_{n-1,\tau}\mathrm{d}\tau\\
&-\mathsf{I}\{Z(\tau-)=n\}[\lambda_{n,\tau}+\mu_{n,\tau}]\mathrm{d}\tau\\
&+\mathsf{I}\{Z(\tau-)=n+1\}\mathrm{d}\boldsymbol{M}_{\mathrm{M}_{n+1}(\tau)}\\
&+\mathsf{I}\{Z(\tau-)=n-1\}\mathrm{d}\boldsymbol{M}_{\Lambda_{n-1}(\tau)}\\
&-\mathsf{I}\{Z(\tau-)=n\}\mathrm{d}[\boldsymbol{M}_{\Lambda_{n}(\tau)}+\boldsymbol{M}_{\mathrm{M}_{n}(\tau)}],\\
&\quad n=1,2,\ldots
\end{aligned}
\]
and
\[
\begin{aligned}
\mathrm{d}\mathsf{I}\{Z(\tau)=0\}=&\mathsf{I}\{Z(\tau-)=1\}\mu_{1,\tau}\mathrm{d}\tau-\mathsf{I}\{Z(\tau-)=0\}\lambda_{0,\tau}\mathrm{d}\tau\\
&+\mathsf{I}\{Z(\tau-)=1\}\mathrm{d}\boldsymbol{M}_{\mathrm{M}_{1}(\tau)}-\mathsf{I}\{Z(\tau-)=0\}\mathrm{d}\boldsymbol{M}_{\Lambda_{0}(\tau)}.
\end{aligned}
\]

Now we rewrite the last two equations in their integral forms by taking the expectation and averaging. We have:

 \[
\begin{aligned}
0=&\lim_{T\to\infty}\frac{1}{T}\mathsf{E}\int_{0}^{T}\mathsf{I}\{Z(\tau-)=n+1\}\mu_{n+1,\tau}\mathrm{d}\tau\\
&+\lim_{T\to\infty}\frac{1}{T}\mathsf{E}\int_{0}^{T}\mathsf{I}\{Z(\tau-)=n-1\}\lambda_{n-1,\tau}\mathrm{d}\tau\\
&-\lim_{T\to\infty}\frac{1}{T}\mathsf{E}\int_{0}^{T}\mathsf{I}\{Z(\tau-)=n\}[\lambda_{n,\tau}+\mu_{n,\tau}]\mathrm{d}\tau\\
&+\lim_{T\to\infty}\frac{1}{T}\mathsf{E}\int_{0}^{T}\mathsf{I}\{Z(\tau-)=n+1\}\mathrm{d}\boldsymbol{M}_{\mathrm{M}_{n+1}(\tau)}\\
&+\lim_{T\to\infty}\frac{1}{T}\mathsf{E}\int_{0}^{T}\mathsf{I}\{Z(\tau-)=n-1\}\mathrm{d}\boldsymbol{M}_{\Lambda_{n-1}(\tau)}\\
&-\lim_{T\to\infty}\frac{1}{T}\mathsf{E}\int_{0}^{T}\mathsf{I}\{Z(\tau-)=n\}\mathrm{d}[\boldsymbol{M}_{\Lambda_{n}(\tau)}+\boldsymbol{M}_{\mathrm{M}_{n}(\tau)}],\\
&\quad n=1,2,\ldots
\end{aligned}
\]
and
\[
\begin{aligned}
0=&\lim_{T\to\infty}\frac{1}{T}\mathsf{E}\int_{0}^{T}\mathsf{I}\{Z(\tau-)=1\}\mu_{1,\tau}\mathrm{d}\tau\\
&-\lim_{T\to\infty}\frac{1}{T}\mathsf{E}\int_{0}^{T}\mathsf{I}\{Z(\tau-)=0\}\lambda_{0,\tau}\mathrm{d}\tau\\
&+\lim_{T\to\infty}\frac{1}{T}\mathsf{E}\int_{0}^{T}\mathsf{I}\{Z(\tau-)=1\}\mathrm{d}\boldsymbol{M}_{\mathrm{M}_{1}(\tau)}\\
&-\lim_{T\to\infty}\frac{1}{T}\mathsf{E}\int_{0}^{T}\mathsf{I}\{Z(\tau-)=0\}\mathrm{d}\boldsymbol{M}_{\Lambda_{0}(\tau)}.
\end{aligned}
\]
In these equations, all the terms containing expectation over martingales (such as $\lim_{T\to\infty}T^{-1}\mathsf{E}\int_{0}^{T}\mathsf{I}\{Z(\tau-)=n-1\}\mathrm{d}\boldsymbol{M}_{\Lambda_{n-1}(\tau)}$, for example) are equal to zero, and
from these equations we finally arrive at the system of equations:
\[
0=P_{n+1}\mu_{n+1}+P_{n-1}\lambda_{n-1}-P_n(\lambda_n+\mu_n), \quad n=1,2,\ldots
\]
and
\[
0=P_1\mu_1-P_0\lambda_0,
\]
where $P_n=\lim_{\tau\to\infty}\mathsf{P}\{Z(\tau)=n\}$ are the final probabilities, if they exist. It is readily seen that the system of equations for the final probabilities coincides with the standard system of equations for the ordinary birth-and-death process.

\subsection{Criteria for recurrence and transience}\label{S2.2}

In the proof of the theorem we use a simplified (degenerate) version of the criteria for recurrence or transience of birth-and-death processes \cite{A1}, which is sufficient for the purpose of this paper. It is as follows.

\begin{lem}\label{6}
An ordinary birth-and-death process is transient, if there exists $c>1$ and $n_0$ such that for all $n\geq n_0$,
\[
\frac{\lambda_{n}}{\mu_{n}}\geq1+\frac{c}{n},
\]
and it is recurrent, if there exists $n_0$ such that for all $n>n_0$
\[
\frac{\lambda_{n}}{\mu_{n}}\leq1+\frac{1}{n}.
\]
\end{lem}
The statement of Lemma \ref{6} follows immediately from the well-known Karlin--McGregor criteria \cite{KM}, according to which an ordinary birth-and-death process is recurrent if and only if $$\sum_{n=1}^{\infty}\prod_{k=1}^{n}\frac{\mu_k}{\lambda_k}=\infty,$$ by application of Raabe's convergence (divergence) test for positive series. The more general results formulated in \cite{A1, A5} were a consequence of the extended version of the Bertrand--De Morgan test \cite{A1} and Bertrand--De Morgan--Cauchy test \cite{A5}.

For the time-inhomogeneous birth-and-death process, the parameters of which obey \eqref{5} -- \eqref{14}, the condition are used to be similar. For transience, it is an existence of $c>1$ and $n_0$ such that for all $n\geq n_0$,
\[
\frac{\lambda_{n,\tau}}{\mu_{n,\tau}}\geq1+\frac{c}{n}.
\]
For recurrence, it is an existence of $c<1$ and $n_0$ such that for all $n\geq n_0$,
\[
\frac{\lambda_{n,\tau}}{\mu_{n,\tau}}\leq1+\frac{c}{n}.
\]
Here we keep in mind that $\tau\to\infty$ together with $n\to\infty$ with probability $1$.

\subsection{The final part of the proof}\label{S2.3} To finish the proof we are aimed to adapt the criteria for recurrence and transience for the time-inhomogeneous birth-and-death process with birth rates $\lambda_{n,\tau}=1/2+\varphi(n,\tau)$ and death rates $\mu_{n,\tau}=1/2-\varphi(n,\tau)$. The birth-and-death process is a continuous time process, while the random walk $X_t$ we deal with in the formulation of the theorem is a discrete time process. To correctly specify our consideration, we assume that we deal with the continuous c\`adl\`ag process $Z(\tau)$, where the moments of jumps characterize the random walk $X_t$. The meaning of $t$ is the $t$th consecutive event of the Poisson processes, the parameters of which are specified at the moments of jumps and depend on state $n$ and discrete time $t$. Let $u_t$ denote the time moment of the $t$th jump. It is assumed that for any $\tau: u_t\leq\tau<u_{t+1}$ the function $\varphi(n,\tau)=\varphi(n,u_t)$, and consequently for the rates of the Poisson processes $\lambda_{n,\tau}=\lambda_{n,u_t}$ and $\mu_{n,\tau}=\mu_{n,u_t}$. Notice that the rates of Poisson processes are specified such that the mean time between two jumps is equal to $1$. It is worth noting that the replacement of the original random walk that is a discrete time process with continuous time-inhomogeneous birth-and-death process does not change the basic property of these two stochastic processes: both of them are either recurrent or transient. The idea of such replacement is not new (e.g. see \cite{A4, A7}).

Denote by $\mathcal{F}_\tau$ the filtration of the process $Z(\tau)$. If $\varphi(Z(\tau),\tau)$ is a positive process vanishing in $L^1$ as $\tau\to\infty$, then denoting by $\triangle Z_{\tau,\tau+\sigma}$ the increment of the process $Z(\tau)$ in the interval $[\tau, \tau+\sigma)$, for the process $Z^2(\tau)-\tau$ we have
\[
\begin{aligned}
\mathsf{E}\left[Z^2(\tau+\sigma)-(\tau+\sigma)~|~\mathcal{F}_{\tau-}\right]&=\mathsf{E}\left[(Z(\tau-)+\triangle Z_{\tau, \tau+\sigma})^2-(\tau+\sigma)~|~\mathcal{F}_{\tau-}\right]\\
&=Z^2(\tau-)+2\mathsf{E}(Z(\tau-)\triangle Z_{\tau, \tau+\sigma}~|~\mathcal{F}_{\tau-})-\tau\\
&=Z^2(\tau)-\tau+o(Z(\tau)),
\end{aligned}
\]
since $\mathsf{E}(Z(\tau-)\triangle Z_{t,t+s}~|~\mathcal{F}_{\tau-})=o(Z(\tau))$, and $\mathsf{E}(\triangle Z_{\tau,\tau+\sigma}^2|\mathcal{F}_{\tau-})=\mathsf{E}\triangle Z_{\tau,\tau+\sigma}^2=\sigma$ for any $\tau$. The last equality follows from Wald's identity \cite{F}. Specifically, we have $\mathsf{E}\triangle Z_{\tau,\tau+\sigma}^2=\mathsf{E}\sum_{i=1}^{N_\sigma}1=\sigma$, where $N_\sigma$ denotes the number of events in time $\sigma$ of a Poisson process with rate $1$.

The derived asymptotic relationship implies that $\varphi(Z(\tau),\tau)\asymp\varphi(Z(\tau),Z^2(\tau))$ for large $\tau$ with probability approaching to1. It also implies that
$\varphi(X_t,t)\asymp\varphi(X_t,X_t^2)$ for large $t$ with probability approaching to $1$.

Then for large $n$ we have
\[
\frac{\lambda_{n,t}}{\mu_{n,t}}\asymp\frac{\lambda_{n,n^2}}{\mu_{n,n^2}}=\frac{1+2\varphi(n,n^2)}{1-2\varphi(n,n^2)}=1+4\varphi(n,n^2)+O[(\varphi(n,n^2))^2].
\]
Now, application of Lemma \ref{6} yields the required statement of the theorem.

\section{Examples}\label{S3}
In this section we provide a number of examples to illustrate the main result of this paper.
\medskip

1. Let $\varphi(n,t)\asymp\rho n/(2t)$, $n\to\infty$ (see \cite{MV}). We show that the process recurrent if $\rho<1/2$, and it is transient if $\rho>1/2$.

Indeed, in this case $\varphi(n, n^2)=\rho/(2n)$ for large $n$. Hence the required statement follows from the theorem.

Let us now find the condition, under which $\varphi(X_t,t)\to0$ in $L^1$. We have
\[
\mathsf{E}\{X_{t+1}~|~X_t\}\asymp X_t\left(1+\frac{\rho}{t}\right), \quad t\to\infty.
\]
From which it is readily seen that $\mathsf{E}X_t=O(t^\rho)$. Hence $\mathsf{E}X_t/t\to0$ as $t\to\infty$ if and only if $\rho<1$. This means that only under this condition $\varphi(X_t,t)\to0$ in $L^1$.

\medskip
2. Let $\varphi(n,t)\asymp n^\alpha/t^\beta$, $\beta>\alpha$, $\beta>0$, $n\to\infty$ (see \cite{MV}).

Assuming $n$ large, for the condition of recurrence we have the inequality
\[
\frac{n^\alpha}{n^{2\beta}}<\frac{c}{4n}, \quad c<1,
\]
which yields $\alpha<2\beta-1$. Keeping in mind that $\alpha<\beta$, with the combination of the obtained result we obtain $\alpha<\min\{\beta, 2\beta-1\}$. The condition of transience is obtained similarly. It is complementary to the condition of recurrence. Namely, it is $0\leq\beta<1$ and $2\beta-1<\alpha<\beta$.

\medskip
3. Let $\varphi(n,t)\asymp\rho n^\alpha/t^\beta$, $-1\leq\alpha\leq1$, $2\beta-\alpha=1$, $n\to\infty$ (unsolved problem in \cite{MV}). We have $\beta=(1+\alpha)/2$, and the asymptotic expression for $\varphi(n,t)$ takes the form
\[
\varphi(n,t)\asymp\frac{\rho n^\alpha}{t^{(1+\alpha)/2}}, \quad n\to\infty.
\]

For large $n$ we obtain
\[
\varphi(n,n^2)\asymp\frac{\rho n^\alpha}{n^{1+\alpha}}=\frac{\rho}{n}.
\]
This enables us to conclude that for $-1\leq\alpha\leq1$ the process $X_t$ is recurrent for $\rho<1/4$ and transient for $\rho>1/4$.

Let us now discuss the condition, under which $\varphi(X_t,t)\to0$ in $L^1$. Here we need only consider the case of $\alpha=1$, since for $-1\leq\alpha<1$ we obviously have the required convergence. For $\alpha=1$ the problem reduces to the case considered in Example 1. Taking into account that result, we arrive at the conclusion that in the case $\alpha=1$, $\varphi(X_t,t)\to0$ in $L^1$ if and only if $\rho<1/2$.

\medskip
4. Let $\varphi(n,t)\asymp\mathrm{e}^{\alpha n-\beta t}$, $\beta>0$, $n\to\infty$.

For large $n$, $\varphi(n,n^2)\asymp\mathrm{e}^{\alpha n-\beta n^2}=o(1/n)$. Hence, the process $X_t$ is recurrent.

\section{Concluding remarks}\label{S4}

In this paper we formulated and proved a new result for recurrence or transience of time-inhomogeneous birth-and-death processes. Compared to the earlier considerations in \cite{MV} our achievements are as follows.

\begin{enumerate}
\item[$\bullet$] We established simple and general criteria for recurrence or transience of time-inhomogeneous birth-and-death processes. This enables us to establish recurrence or transience of the processes in a very simple way for a wide class of processes. Our approach makes an essential difference with \cite{MV}, where special routine derivations for any particular case study were required.
\item[$\bullet$] The generality of our main result enabled us to solve the open problem that was left unsolved in \cite{MV}.
\item[$\bullet$] In Examples 1 and 3 we also studied the behaviour of the drift. We found the value of the parameter, under which the drift was vanishing in time in the sense provided in the paper.
\end{enumerate}

\subsection*{Acknowledgement} The author thanks all the people who helped in preparation of this paper officially or privately.

\end{document}